\documentclass[12pt]{article}
\usepackage{color}
\definecolor{darkblue}{rgb}{0.00,0.25,0.50}
\usepackage[colorlinks,filecolor=blue,citecolor=darkblue]{hyperref}

\usepackage{amsfonts,amssymb,amsmath,amsthm}
\usepackage{url}
\usepackage{enumerate}
\usepackage[ukrainian, russian, english]{babel}
\usepackage[cp1251]{inputenc}
\begin{document} \selectlanguage{ukrainian}
\thispagestyle{empty}

\title{}

\begin{center}
\textbf{\Large Approximative characteristics of classes \\ of functions  \boldmath{$S^{\Omega}_{p,\theta}B(\mathbb{R}^d)$} with a given majorant \\ of mixed modulus of smoothness}
\end{center}

\vskip 3mm

\begin{center}
\textbf{\Large  Апроксимативні характеристики функцій \\ з класів \boldmath{$S^{\Omega}_{p,\theta}B(\mathbb{R}^d)$} із заданою мажорантою  \\ мішаних модулів неперервності}
\end{center}

\vskip0.5cm

\begin{center}
S.~Ya.~Yanchenko,  \ S.~A.~Stasyuk\\ \emph{\small
Institute of Mathematics of NAS of
Ukraine, Kyiv\\}
\end{center}
\begin{center}
C.~Я.~Янченко, \ С.~А.~Стасюк \\
\emph{\small Інститут математики НАН України, Київ}
\end{center}
\vskip0.5cm

\vskip 4mm

\begin{abstract}

 We obtain order estimates of approximation of functions from the classes
${S^{\Omega}_{p,\theta}B (\mathbb{R}^d)}$ in the space $L_q(\mathbb{R}^d)$, $1<p<q<\infty$, by entire functions of
exponential type with supports of their Fourier transforms in sets generated by the
level surfaces of a function $\Omega$.

Одержано точні за порядком оцінки наближення функцій  з класів
$S^{\Omega}_{p,\theta}B (\mathbb{R}^d)$ у просторі $L_q(\mathbb{R}^d)$, $1<p<q<\infty$, за допомогою цілих функцій експоненціального типу з носіями їх перетворень Фур'є на множинах, які породжуються
поверхнями рівня функції $\Omega$.
\end{abstract}

\vskip 0.5 cm


\noindent \textbf{1. Означення класів функцій та апроксимативних характеристик.} У статті продовжено дослідження  апроксимативних характеристик функцій
з  класів Нікольського--Бєсова
$S^{\Omega}_ {p,\theta}B(\mathbb{R}^d)$ у просторі
$L_q(\mathbb{R}^d)$~\cite{Stasuk_Yanchenko_Anal_math}, \cite{Yanchenko_YMG_2016_12}. Встановлено точні за порядком оцінки наближення таких функцій
за допомогою цілих функцій експоненціального типу з
носіями їх перетворень Фур'є на множинах, які породжуються
поверхнями рівня функції $\Omega$ у випадку, коли $1<p<q<\infty$.

Нехай $\mathbb{R}^d$
--- $d$-вимірний евклідів простір з елементами
${\boldsymbol{x}=(x_1,...,x_d)}$ і
${(\boldsymbol{x},\boldsymbol{y})=x_1y_1+...+x_dy_d}$. Через
$L_q(\mathbb{R}^d)$, $1\leqslant  q\leqslant \infty$, позначимо
простір вимірних на $\mathbb{R}^d$ функцій зі
скінченною нормою
 $$
\|f\|_{L_{q}}:=\|f\|_q:=
\left(\int\limits_{\mathbb{R}^{d}}|f(\boldsymbol{x})|^{q}d\boldsymbol{x}
\right) ^{\frac{1}{q}}, \ 1\leqslant  q<\infty,
 $$
 $$
 \|f\|_{L_{\infty}}:= \|f\|_{\infty}:=
 \mathop {\rm ess \sup}\limits_{\boldsymbol{x}\in \mathbb{R}^d}   |f(\boldsymbol{x})|.
 $$

Для функції $f\in L_q{(\mathbb{R}^d)}$ розглянемо різницю $l$-го
порядку,  $l\in\mathbb{N}$, за змінною $x_j$ з кроком $h_j$, яка
визначається таким чином:
$$
\Delta^{l}_{h_{j}}f(\boldsymbol{x}):=\sum
\limits_{n=0}^{l}(-1)^{l-n}C^{n}_{l}f(x_1,...,
x_{j-1},x_j+nh_j,x_{j+1},...,x_d).
$$
Також означимо мішану кратну різницю $l$-го порядку  функції $f$ з векторним кроком
\mbox{$\boldsymbol{h}=(h_1,\dots,h_d)$}:
$$
\Delta^{l}_{\boldsymbol{h}}f(\boldsymbol{x})=
\Delta^{l}_{h_d}\big(\Delta^{l}_{h_{d-1}}...(\Delta^{l}_{h_{1}}f(\boldsymbol{x}))\big).
$$

Мішаний модуль неперервності порядку $l$ функції $f \in L_q(\mathbb{R}^d)$ визначається згідно з формулою
$$
\Omega_{l}(f,\boldsymbol{t})_{q}:=\sup \limits_{|\boldsymbol{h}|\leqslant  \boldsymbol{t}} \|\Delta^{l}_{\boldsymbol{h}}f(\cdot)\|_{q},
$$
де $|\boldsymbol{h}|=(|h_1|,...,|h_d|)$, а нерівності типу $\boldsymbol{a}\leqslant \boldsymbol{b}$ $(\boldsymbol{a}>\boldsymbol{b})$ для векторів ${\boldsymbol{a}=(a_1,...,a_d)}$ та ${\boldsymbol{b}=(b_1,...,b_d)}$ тут і надалі розуміємо покоординатно, тобто $a_j\leqslant b_j$ $(a_j>b_j)$, $j=\overline{1,d}$. Також будемо використовувати запис $\boldsymbol{t}\geqslant \boldsymbol{0}$, якщо ${t_j\geqslant{0}}$,
$j=\overline{1,d}$.

Нехай $\Omega(\boldsymbol{t})$, $\boldsymbol{t}=(t_1,\dots,t_d)$~--- функція
типу мiшаного модуля неперервності порядку $l$, тобто функція, яка визначена і неперервна
 на $\mathbb{R}^d_{+}$, що задовольняє такі умови:
\begin{itemize}
 \item[1)] $\Omega(\boldsymbol{t})>0$, $\boldsymbol{t}>\boldsymbol{0}$ і
$\Omega(\boldsymbol{t})=0$, якщо $\prod \limits_{j=1}^dt_j=0$;

\item[2)] $\Omega(\boldsymbol{t})$ неспадна за кожною змiнною;

\item[3)] $\Omega(m_1t_1,\dots,m_dt_d)\leqslant  \left(\prod \limits_{j=1}^dm_j\right)^l
\Omega(\boldsymbol{t})$, $m_j\in \mathbb{N},$ $j=\overline{1,d}$.

\end{itemize}

Множину таких функцій $\Omega$ позначимо через
$\Psi_{l}$.

Додатково будемо вимагати, щоб функція $\Omega$
задовольняла умови $(S^{\alpha})$ та $(S_l)$, які називають умовами
Барi--Стєчкiна \cite{Bari_Stechkin}. Сформулюємо їх:
\begin{itemize}
\item[а)] функція однієї змінної $\varphi(\tau)\geqslant 0$ задовольняє
умову $(S^{\alpha})$, якщо існує таке $\alpha>0$, що $\varphi(\tau)/\tau^{\alpha}$  майже
зростає, тобто існує така незалежна від $\tau_1$  i
$\tau_2 $ стала ${C_1>0}$, що
 $$
 \frac{\varphi(\tau_1)}{\tau_1^{\alpha}} \leqslant C_1\frac{\varphi(\tau_2)}{\tau_2^{\alpha}},
  \ \ \ \ 0<\tau_1 \leqslant \tau_2 \leqslant 1;
 $$

\item[б)] функція однієї змінної $\varphi(\tau)\geqslant 0$ задовольняє
умову $(S_l)$, якщо  існує таке $\gamma$, $0<\gamma<l$, що
$\varphi(\tau)/\tau^{l-\gamma}$ майже спадає, тобто
існує така незалежна від $\tau_1$ i $\tau_2 $ стала $C_2>0$,  що
$$
\frac{\varphi(\tau_1)}{\tau_1^{l-\gamma}} \geqslant
C_2\frac{\varphi(\tau_2)}{\tau_2^{l-\gamma}},  \ \ \ \ 0<\tau_1
\leqslant \tau_2 \leqslant 1 .
$$
\end{itemize}

Будемо вважати, що $\Omega$ задовольняє умови $(S^{\alpha})$
та $(S_l)$, якщо $\Omega$ задовольняє ці умови за
кожною змінною \ $t_j$ при всіх фіксованих значеннях змінних $t_i$, \ ${i\neq j}$. У
випадку, коли для $\Omega$ виконана умова
$(S^{\alpha})$, будемо говорити, що $\Omega$ належить множині
$S^{\alpha}$, а якщо умова $(S_l)$, то~--- множині $S_l$. Стверджуючи це
(також і для функції $\omega$ однієї змінної),
використовуватимемо запис $\Omega \in \Phi_{\alpha,
l}$, ($\omega \in \Phi_{\alpha, l}$), $l\in \mathbb{N}$, де
множина $\Phi_{\alpha, l}$ визначається співвідношенням
\mbox{$\Phi_{\alpha, l}=\Psi_l\cap S^{\alpha} \cap S_l$}.

Зазначимо, що до множини $\Phi_{\alpha, l}$ належать,
наприклад, функції
$$
\Omega(\textbf{t})=\Omega(t_1,\dots,t_d)=\begin{cases}     \prod
\limits_{j=1}^{d}\frac{t_j^{r_j}}{\left\{\log \frac
                                {1}{t_j}\right\}_{+}^{b_j}}, &  \mbox{при} \ \ t_j>0, \ j=\overline{1,d}; \\
    0, & \mbox{при} \ \ \prod \limits_{j=1}^{d}t_j=0,
 \end{cases}
$$
де $\left\{\log \tau \right\}_{+}=\max\left\{ 1; \ \log_2 \tau
\right\}$,  $r_j, b_j\in \mathbb{R}$, $0<r_j<l$, $j=\overline{1,
d}$.

Нехай, далі, $e_d:=\{1,2, ... ,d\}$,  ${d\in \mathbb{N}}$, і
${e:=\{j_1, ... ,j_m \}}$, ${m\leqslant  d}$, ${m\in \mathbb{N}}$, \
${1\leqslant  j_1<j_2< ... <j_m\leqslant d}$,
$\boldsymbol{t}^e=(t_{j_1},\dots,t_{j_m})$,
${\bar{\boldsymbol{t}}^e:=(\bar{t}_1,\dots,\bar{t}_d)}$, де
$$
 \bar{t}_i=
 \begin{cases}
    t_i, & i \in e, \\
    1, & i \in e_d\backslash e.
 \end{cases}
$$

Простори $S^{\Omega}_{p,\theta}B(\mathbb{R}^d)$ для   $1\leqslant  p, \theta
\leqslant  \infty$  і
$\Omega \in \Psi_{l}$ означаються таким чином (див.,
наприклад, \cite{Stasuk_Yanchenko_Anal_math})
$$
S^{\Omega}_{p,\theta}B(\mathbb{R}^d):=\Big\{f\in L_p(\mathbb{R}^d):
\ \|f\|_{S^{\Omega}_{p,\theta}B(\mathbb{R}^d)}< \infty \Big\},
$$
де
 $$
  \|f\|_{S^{\Omega}_{p,\theta}B(\mathbb{R}^d)}:=\|f\|_p+
  \sum \limits_{ e \subset e_d \atop e \neq \varnothing} \left(\int \limits_{0}^{2}\cdot\cdot\cdot\int \limits_{0}^{2}
  \left(\frac {\Omega_{l^e}(f,\boldsymbol{t}^e)_p}{\Omega(\bar{\boldsymbol{t}}^e)}\right)^{\theta}
   \prod \limits_{j\in e}\frac{dt_j}{t_j}\right)^{\frac{1}{\theta}},
$$
якщо $1\leqslant  \theta  <\infty$, та
$$%
  \|f\|_{S^{\Omega}_{p,\infty}B(\mathbb{R}^d)}:= \|f\|_p+\sum \limits_{e \subset e_d \atop e \neq \varnothing}
  \sup \limits_{\boldsymbol{t}^e>\boldsymbol{0}}   \frac {\Omega_{l^e}(f,\boldsymbol{t}^e)_p}{\Omega(\bar{\boldsymbol{t}}^e)},
$$
де
$$
\Omega_{\boldsymbol{l}^e}(f,\boldsymbol{t}^e)_{q}:= \sup
\limits_{|\boldsymbol{h}^e|\leqslant \boldsymbol{t}^e}
\|\Delta^{\boldsymbol{l}^e}_{\boldsymbol{h}^e}f(\boldsymbol{x})\|_{q}, \ \
\boldsymbol{h}^e:=(h_{j_1},\dots,h_{j_m}),
$$
$$
\Delta^{\boldsymbol{l}^e}_{\boldsymbol{h}^e}f(\boldsymbol{x})=
\Delta^{l}_{h_{j_m}}\big(\Delta^{l}_{h_{j_{m-1}}}...
(\Delta^{l}_{h_{j_1}}f(\dots,x_{j_1},\dots,x_{j_m},\dots))\big).
$$

Зауважимо, що простори функцій
$S^{\Omega}_{p,\theta}B(\mathbb{R}^d)$  є узагальненням відомих
просторів $S^{\boldsymbol{r}}_{p,\theta}B(\mathbb{R}^d)$~\cite{Nikolsky_63}, \cite{Amanov_1965}, що визначаються при явному заданні функції $\Omega$, а саме ${\Omega(\boldsymbol{t})=\boldsymbol{t}^{\boldsymbol{r}}=t_1^{r_1}\cdot\ldots\cdot
t_d^{r_d}}$, $0<r_j<l$, $j=\overline{1, d}$.
Нагадаємо, що простори ${S^{\boldsymbol{r}}_{p}H(\mathbb{R}^d)=}$ $={S^{\boldsymbol{r}}_{p,\infty}B(\mathbb{R}^d)}$
були вперше розглянуті С.\,М.~Нікольським~\cite{Nikolsky_63}, простори $S^{\boldsymbol{r}}_{p,\theta}B(\mathbb{R}^d)$, при $1\leqslant \theta <\infty$ були введені
Т.\,І.~Амановим~\cite{Amanov_1965} (див. також \cite{Amanov_book}).
Надалі  будемо використовувати скорочені позначення
$S^{\Omega}_{p,\theta}B$, $S^{\boldsymbol{r}}_{p,\theta}B$ і
$S^{\boldsymbol{r}}_{p}H$ відповідно для
$S^{\Omega}_{p,\theta}B(\mathbb{R}^d)$, $S^{\boldsymbol{r}}_{p,\theta}B(\mathbb{R}^d)$ та
$S^{\boldsymbol{r}}_{p}H(\mathbb{R}^d)$ .

Дослідження класів Нікольського--Бєсова з домінуючою мішаною похідною   $S^{\boldsymbol{r}}_{p,\theta}B$, з точки зору знаходження порядкових оцінок деяких апроксимативних характеристик проводилися, зокрема, у роботах  Wang Heping і Sun Yongsheng~\cite{WangHeping_SunYongsheng_1995}, Wang Heping~\cite{WangHeping_2004}. З основними результатами  щодо дослідження класів Нікольського--Бєсова з домінуючою мішаною похідною  у періодичному випадку можна ознайомитися в монографії В.\,Н.~Темлякова~\cite{Temlyakov_1986m}, якщо $\theta=\infty$ (для класів Нікольського), та у монографії А.\,С.~Романюка~\cite{Romanyuk_2012m}, якщо $1\leqslant \theta < \infty$ (для класів Бєсова). На даний час є значний інтерес до дослідження різних аналогів класів Нікольського--Бєсова, які визначаються гладкісним параметром $\Omega$, що підпорядкований  деяким додатковим умовам: Н.\,Н.\,Пустовойтов~\cite{Pustovoitov_94}, \cite{Pustovoitov_2012}, Wang Heping і Sun Yongsheng~\cite{WangHeping_SunYongsheng_1997},  Liqin Duan~\cite{Liqin_Duan2010} та ін.

В \cite{Stasuk_Yanchenko_Anal_math} встановлено еквівалентне нормування лінійних просторів $S^{\Omega}_{p,\theta}B$ опосередковано через, так зване, декомпозиційне представлення елементів цих просторів (див. нижче Теорему~A).  Зазначимо, що  для просторів Нікольського--Бєсова функцій мішаної гладкості, вперше декомпозиційне представлення та відповідне йому нормування з'явилося у роботі   С.\,М.~Нікольського та П.\,І.~Лізоркіна~\cite{Lizorkin_Nikolsky_1989} і, як з'ясувалося пізніше, відіграло ключову роль у дослідженнях, які пов'язані з апроксимацією класів функцій. Оскільки у формулюванні результату з~\cite{Stasuk_Yanchenko_Anal_math}, щодо нормування простору  $S^{\Omega}_{p,\theta}B$, присутні величини, які означаються за допомогою перетворення Фур'є функцій, що визначені на $\mathbb{R}^d$, то наведемо відповідні означення.

Нехай $S=S(\mathbb{R}^d)$~--- простір Л.~Шварца основних нескінченно
диференційовних на $\mathbb{R}^d$ комплекснозначних функцій
$\varphi$, що спадають на нескінченності разом зі своїми похідними
швидше за будь-який степінь функції $\left(x_1^2+\ldots+x_d^2\right)^{-\frac{1}{2}}$ (див., наприклад,
\cite{Lizorkin_69}). Через $S'$ позначимо
простір лінійних неперервних функціоналів над $S$. Зазначимо, що
елементами простору $S'$ є узагальнені функції повільного росту. Якщо $f\in S'$, то $\langle f,\varphi\rangle$ позначає значення функціонала $f$ на пробній функції
$\varphi \in S$.

Перетворення Фур'є $\mathfrak{F}\varphi: S\rightarrow S$
визначається згідно з формулою:
$$
(\mathfrak{F}\varphi)(\boldsymbol{\lambda})=\frac{1}{(2\pi)^{d/2}}\int
 \limits_{\mathbb{R}^d}\varphi(\boldsymbol{t})
 e^{-i(\boldsymbol{\lambda},\boldsymbol{t})}d\boldsymbol{t}
 \equiv \widetilde{\varphi}(\boldsymbol{\lambda}).
$$

Обернене перетворення Фур'є задається таким чином:
$$
(\mathfrak{F}^{-1}\varphi)(\boldsymbol{t})=\frac{1}{(2\pi)^{d/2}}
\int \limits_{\mathbb{R}^d}\varphi(\boldsymbol{\lambda})
e^{i(\boldsymbol{\lambda},\boldsymbol{t})}d\boldsymbol{\lambda}\equiv
\widehat{\varphi}(\boldsymbol{t}).
$$

Перетворення Фур'є (обернене перетворення  Фур'є) узагальнених функцій $f\in S'$  визначається згідно з формулою
$$
\langle \mathfrak{F}f,\varphi\rangle=\langle f,\mathfrak{F}\varphi
\rangle,  \ \ \ \langle \widetilde{f},\varphi\rangle=\langle
f,\widetilde{\varphi} \rangle, \ \ \varphi \in S,
$$
$$
(\langle \mathfrak{F}^{-1}f,\varphi\rangle=\langle
f,\mathfrak{F}^{-1}\varphi \rangle,  \ \ \ \langle
\widehat{f},\varphi\rangle=\langle f,\widehat{\varphi} \rangle, \ \ \varphi \in S).
$$

Носієм узагальненої функції $f$ будемо називати замикання
$\overline{\mathfrak{N}}$ такої множини точок
$\mathfrak{N}\subset\mathbb{R}^d$, що для довільної $\varphi \in S$,
яка дорівнює нулю в $\overline{\mathfrak{N}}$, виконується рівність
$\langle f,\varphi \rangle = 0$. Носій узагальненої функції $f$
будемо позначати через $\mbox{supp}\, f$. Також будемо говорити, що функція $f$
 зосереджена на множині $G$, якщо $\mbox{supp}\, f \subseteq G$.

Зазначимо, що для $1\leqslant p \leqslant \infty$ існує природне неперервне
вкладення $L_p(\mathbb{R}^d)$ в $S'$ і в цьому сенсі функції з
$L_p(\mathbb{R}^d)$ ототожнюються з елементами з $S'$.

Далі, для кожного вектора ${\boldsymbol{s}=(s_{1},...,s_{d})}$,
${s_{j}\in \mathbb{Z}_+,  \ j=\overline{1,d}}$, розглянемо множину
$$
Q_{2^{\boldsymbol{s}}}^*\!:=Q^*(\boldsymbol{s})\!:=
\Big\{\boldsymbol{\lambda}=(\lambda_1,\ldots,\lambda_d): \
\eta(s_j)2^{s_{j}-1}\leqslant   |\lambda_j|<2^{s_j},
$$
$$ \lambda_j
\in \mathbb{R}, \ j=\overline{1,d}\Big\},
$$
де $\eta(0)=0$ і $\eta(t)=1$, $t>0$.

Нехай $\mathcal{A}\subset\mathbb{R}^d$~--- деяка вимірна множина. Позначимо
через $\chi_\mathcal{A}$ характеристичну функцію множини
$\mathcal{A}$ і для $f\in L_p(\mathbb{R}^d)$ покладемо
$$
\delta_{\boldsymbol{s}}^*(f,\boldsymbol{x})=\mathfrak{F}^{-1}\big(\chi_{Q_{2^{\boldsymbol{s}}}^*}
\cdot
\mathfrak{F}f\big).
$$

\bf Теорема A \rm (\cite{Stasuk_Yanchenko_Anal_math}). \it Нехай $1< p < \infty$ і $\Omega
\in \Phi_{\alpha, l}$. Функція $f$ належить простору
$S^{\Omega}_{p,\theta}B$,  $1\leqslant\theta<\infty$, тоді і тільки
тоді, коли
 $$
  \left\{\sum\limits_{\boldsymbol{s}\geqslant 0}
    \big(\Omega(2^{-{\boldsymbol{s}}})\big)^{-\theta}\|\delta^*_{\boldsymbol{s}}(f,\cdot)\|_p^\theta\
   \right\}^{\frac{1}{\theta}}<\infty,
 $$
до того ж
 \begin{equation}\label{norm_1_dek}
   \|f\|_{S^{\Omega}_{p,\theta}B}\asymp   \left\{\sum\limits_{\boldsymbol{s}\geqslant 0}
   \big(\Omega(2^{-\boldsymbol{s}})\big)^{-\theta}\|\delta^*_{\boldsymbol{s}}(f,\cdot)\|_p^\theta\ \right\}^{\frac{1}{\theta}},
 \end{equation}
 де $\Omega(2^{-\boldsymbol{s}})=\Omega(2^{-s_1},\ldots,2^{-s_d})$.

Функція $f$ належить простору $S^{\Omega}_{p,\infty}B$, тоді і
тільки тоді, коли
 $$
  \sup\limits_{\boldsymbol{s}\geqslant 0} \frac{\|\delta^*_{\boldsymbol{s}}(f,\cdot)\|_p}
  {\Omega(2^{-\boldsymbol{s}})}   < \infty,
 $$
до того ж
 \begin{equation}\label{norm_inf_dek}
  \|f\|_{S^{\Omega}_{p,\infty}B}\asymp
  \sup\limits_{\boldsymbol{s}\geqslant0} \frac{\|\delta^*_{\boldsymbol{s}}(f,\cdot)\|_p}{\Omega(2^{-\boldsymbol{s}})} .
 \end{equation}
\rm

 Тут і надалі по тексту для додатних величин $A$ і  $B$ використовується запис  $A\asymp B$, який означає, що існують такі додатні  сталі $C_3$ та $C_4$, які не залежать від одного істотного параметра у величинах  $A$ і  $B$ (наприклад, у вище наведених співвідношеннях (\ref{norm_1_dek}) і (\ref{norm_inf_dek})~--- від функції $f$), що ${C_3 A \leqslant B \leqslant C_4 A}$. Якщо тільки $B\leqslant C_4 A $ $\big(B \geqslant C_3 A\big)$, то пишемо
$B\ll A$ $\big(B \gg A \big)$. Всі сталі $C_i$, $i=1,2,...$, які
зустрічаються у роботі, залежать, можливо, лише від параметрів, що
входять в означення класу, метрики, в якій оцінюється похибка
наближення, та розмірності простору $\mathbb{R}^d$.

Під класом $S^{\Omega}_{p,\theta}B$ будемо розуміти множину
функцій ${f \in L_p(\mathbb{R}^d)}$ для яких
$\|f\|_{S^{\Omega}_{p,\theta}B}\leqslant  1$ і при цьому збережемо
для класів $S^{\Omega}_{p,\theta}B$ ті ж самі позначення, що і для
просторів $S^{\Omega}_{p,\theta}B$.

Перейдемо до означення апроксимативних характеристик.

Нехай $\mathcal{L}\subset \mathbb{Z}^d_+$~--- деяка обмежена
множина. Покладемо
$$
Q(\mathcal{L})=\bigcup\limits_{\boldsymbol{s}\in \mathcal{L}}Q^*(\boldsymbol{s})
$$
і позначимо
$$
G\big(Q(\mathcal{L})\big)=\Big\{f\in L_q(\mathbb{R}^d)\colon
\mbox{supp}\mathfrak{F}f\subseteq Q(\mathcal{L}) \Big\}.
$$
Відомо, що елементами множини $G\big(Q(\mathcal{L})\big)$ є цілі функції експоненціального типу.

Для $f\in L_q(\mathbb{R}^d)$, $1\leqslant q \leqslant \infty$,
означимо величину
$$
E\big(f, G\big(Q(\mathcal{L})\big)\big)_q:=E_{Q(\mathcal{L})}(f)_q:=
\inf\limits_{g\in  G(Q(\mathcal{L}))}\|f(\cdot)-g(\cdot)\|_q,
$$
яка називається найкращим наближенням функції $f$ цілими функціями з
множини $G\big(Q(\mathcal{L})\big)$. Якщо $F\subset
L_q(\mathbb{R}^d)$~--- деякий функціональний клас, то покладемо
\begin{equation}\label{EQN}
E_{Q(\mathcal{L})}(F)_q=\sup\limits_{f\in F}E_{Q(\mathcal{L})}(f)_q.
\end{equation}

Далі для $f\in L_q(\mathbb{R}^d)$, $1\leqslant q \leqslant \infty$, покладемо
$$
S_{Q(\mathcal{L})}f(\boldsymbol{x})= S_{Q(\mathcal{L})}(f,\boldsymbol{x})=\sum\limits_{\boldsymbol{s}\in
\mathcal{L}}\delta_{\boldsymbol{s}}^*(f,\boldsymbol{x}), \ \ \boldsymbol{x}\in  \mathbb{R}^d
$$
і означимо
\begin{equation}\label{EQN1}
\mathcal{E}_{Q(\mathcal{L})}(f)_q=
\|f(\cdot)-S_{Q(\mathcal{L})}f(\cdot)\|_q,  \ \  \mathcal{E}_{Q(\mathcal{L})}(F)_q=\sup\limits_{f\in F}\mathcal{E}_{Q(\mathcal{L})}(f)_q.
\end{equation}

Наші дослідження величин (\ref{EQN}) та (\ref{EQN1}) проводяться у випадку, коли ${F=S^{\Omega}_{p,\theta}B}$, а множина $\mathcal{L}$ певним чином пов'язана з функцією $\Omega$.

Для будь-якого $N\in \mathbb{N}$, $1<p<q<\infty$ покладемо
$$
\kappa(\Omega,
N):=\kappa(N):=\left\{\boldsymbol{s}=(s_1,\ldots,s_d)\in\mathbb{Z}^d_+ \colon \Omega(2^{-\boldsymbol{s}})2^{\|\boldsymbol{s}\|_1\left(\frac{1}{p}-\frac{1}{q}\right)}
\geqslant\frac{1}{N}\right\},
$$
$$
Q(\kappa(N)):=Q(N)=:\bigcup\limits_{\boldsymbol{s}\in\kappa(N)}
Q^*(\boldsymbol{s}),
$$
де $\|\boldsymbol{s}\|_1=s_1+\ldots+s_d$.

Зазначимо, що множини $Q(N)$ породжуються поверхнями рівня функції $\Omega(\boldsymbol{t})\boldsymbol{t}^{-\left(\frac{1}{p}-\frac{1}{q}\right)}$, $\Omega(\boldsymbol{t})\in \Phi_{\alpha,l}$, $\alpha>\frac{1}{p}-\frac{1}{q}$. Якщо
$$
\Omega(\boldsymbol{t})=\Omega_1(\boldsymbol{t})/\prod\limits_{j=1}^{d}
t_j^{-\left(\frac{1}{p}-\frac{1}{q}\right)}
$$ і
$\Omega_1(\boldsymbol{t})=\prod\limits_{j=1}^{d}t_j^{r_j}$, $0<r_j<l$, $j=\overline{1,d}$, одержимо множини $Q(N)$, які називаються східчастими гіперболічними хрестами.

Зазначимо, що наближення класів Нікольського--Бєсова періодичних функцій мішаної гладкості тригонометричними поліномами зі спектром у східчастому гіперболічному хресті та у множинах $Q(N)$ розглядалися, зокрема, у роботах \cite{D_Zung_86}, \cite{Romanyuk_91}, \cite{Romanyuk_92}, \cite{Pustovoitov_99}, \cite{Stasyuk_2010mz}, \cite{Stasyuk_2014Tr}, \cite{Balgimbaeva_2015Tr}. З детальною історією питання можна ознайомитися в оглядовій статті~\cite{Cross_2016}.

Формулювання допоміжних результатів та доведення основних потребує означення ще деяких множин в $\mathbb{Z}^d_+$.

Покладемо
$$
\kappa^{\perp}(\Omega, N)\!:=\!\kappa^{\perp}(N)\!:=\!
\left\{\!\boldsymbol{s}=(s_1,\ldots, s_d)\in\mathbb{Z}^d_+\colon
\Omega(2^{-\boldsymbol{s}})2^{\|\boldsymbol{s}\|_1\left(\frac{1}{p}-\frac{1}{q}\right)}\!
<\! \frac{1}{N}\!\right\},
$$
$$
\Theta(N):=\kappa^{\perp}(N) \backslash  \kappa^{\perp}(2^l N),
$$
 тобто
\begin{equation}\label{Omega_1n}
\Theta(N)=\Big\{ \boldsymbol{s}=(s_1,\ldots, s_d)\in\mathbb{Z}^d_+\colon
\frac{1}{2^lN}\leqslant\Omega(2^{-\boldsymbol{s}})2^{\|\boldsymbol{s}\|_1
\left(\frac{1}{p}-\frac{1}{q}\right)}<\frac{1}{N}\Big\}
\end{equation}

У \cite{Pustovoitov_99} показано, що має місце співвідношення
\begin{equation}\label{log}
|\Theta(N)|\asymp \big(\log_2 N\big)^{d-1},
\end{equation}
де через $|A|$ позначено кількість елементів скінченної множини $A$.

Мають місце такі твердження.

 \bf Лема А \rm(\cite{Pustovoitov_99}).  \it Нехай
$\Omega$~--- функція типу мішаного модуля
неперервності порядку $l$, яка задовольняє умову $(S^{\alpha})$, $\alpha>0$. Тоді для
$0<\mu<\infty$
\begin{equation}\label{Omega_lemma}
\sum \limits_{\boldsymbol{s}\in \kappa^{\perp}(N)}(\Omega(2^{-\boldsymbol{s}}))^{\mu} \ll  \sum
\limits_{\boldsymbol{s}\in \Theta(N)}(\Omega(2^{-\boldsymbol{s}}))^{\mu} .
\end{equation}
\rm

\bf Лема Б \rm (\cite{Pustovoitov_99}). \it Нехай $\Omega$~--- функція типу мішаного модуля
неперервності порядку $l$, яка задовольняє умову $(S^{\alpha})$, при $\alpha>\beta>0$. Тоді для
$0<\mu<\infty$
\begin{equation}\label{Omega_pq_lemma}
\sum \limits_{\boldsymbol{s}\in \kappa^{\perp}(N)}\left(\Omega(2^{-\boldsymbol{s}})2^{\|s\|_1\beta}\right)^{\mu} \ll  \sum
\limits_{\boldsymbol{s}\in \Theta(N)}\left(\Omega(2^{-\boldsymbol{s}})2^{\|s\|_1\beta}\right)^{\mu} .
\end{equation}
\rm

Як наслідок з
(\ref{Omega_lemma}), (\ref{Omega_pq_lemma}), (\ref{Omega_1n}) та (\ref{log}) для ${0<\mu<\infty}$ маємо
$$
\sum \limits_{\boldsymbol{s}\in
\kappa^{\perp}(N)}\left(\Omega(2^{-\boldsymbol{s}})2^{\|\boldsymbol{s}\|_1
\left(\frac{1}{p}-\frac{1}{q}\right)}\right)^{\mu} \ll  \sum
\limits_{\boldsymbol{s}\in
\Theta(N)}\left(\Omega(2^{-\boldsymbol{s}})2^{\|\boldsymbol{s}\|_1
\left(\frac{1}{p}-\frac{1}{q}\right)}\right)^{\mu}<
$$
\begin{equation}\label{1n_log}
<\left(\frac{1}{N}\right)^{\mu} \sum \limits_{\boldsymbol{s}\in
\Theta(N)} 1=\left(\frac{1}{N}\right)^{\mu} |\Theta(N)|\asymp \left(\frac{1}{N}\right)^{\mu} \big(\log_2 N\big)^{d-1}.
\end{equation}

\bf Теорема Б \rm (Літтлвуда\,--\,Пелі) (див.,
наприклад, \cite[c.~81]{Nikolsky_1969_book}, \cite{Lizorkin_1967}). \it  Нехай задано
${1<p<\infty}$. Існують такі додатні числа $C_{5}, C_{6}$,  що для
кожної функції ${f\in L_{p}(\mathbb{R}^d)}$ виконуються
співвідношення
 $$
 C_{5}\|f\|_{p}\leqslant  \Bigg\| \Bigg(\sum\limits_{\boldsymbol{s}\geqslant 0}
  |\delta_{\boldsymbol{s}}^*(f,\cdot)|^{2} \Bigg)^{\frac{1}{2}} \Bigg\|_{p}\leqslant  C_{6}\|f\|_{p}.
$$
\rm

У подальших міркуваннях нами буде використана лема, яка є аналогом відповідної леми у періодичному випадку
\cite[гл.~1, \S~3]{Temlyakov_1986m}.

\bf Лема В \rm (\cite{WangHeping_SunYongsheng_1995}). \it Нехай задано $1<p<q<\infty$ і $f\in
L_p(\mathbb{R}^d)$. Тоді
 $$
\|f\|_q\ll \bigg( \sum \limits_{\boldsymbol{s}\geqslant 0} \|\delta_{\boldsymbol{s}}^*(f,\cdot)\|_p^q
\; 2^{\|\boldsymbol{s}\|_1\left(\frac{1}{p}-\frac{1}{q}\right)q}\bigg)^{\frac{1}{q}}.
 $$
\rm

\textbf{2. Порядкові оцінки апроксимативних характеристик класів \boldmath{$S^{\Omega}_{p,\theta}B(\mathbb{R}^d)$} у метриці простору \boldmath{$L_q(\mathbb{R}^d)$, $1<p<q<\infty$}.}
Наша мета полягає у встановленні порядкових по параметру $N$ оцінок величин $E_{Q(\mathcal{L})}(F)_q$ та $\mathcal{E}_{Q(\mathcal{L})}(F)_q$ для $F=S^{\Omega}_{p,\theta}B$ у випадку, коли $\mathcal{L}=\kappa(N)$, при певних обмеженнях на параметри $p$, $q$, $\theta$ та $\Omega$.

При цьому зазначимо, що при $1< q < \infty$ і  $f\in L_q(\mathbb{R}^d)$
 має  місце
співвідношення (див., наприклад,~\cite{Lizorkin_Nikolsky_1989})
\begin{equation}\label{e=e}
E_{Q(\mathcal{L})}(f)_q\leqslant\mathcal{E}_{Q(\mathcal{L})}(f)_q\leqslant C E_{Q(\mathcal{L})}(f)_q,
\end{equation}
де $C\geqslant1$~---  деяка стала.

\bf Теорема 1. \rm \it Нехай $1< p<q< \infty$,
 $1 \leqslant \theta \leqslant \infty$, $\Omega\in \Phi_{\alpha,l}$, з деяким $\alpha> \frac{1}{p}-\frac{1}{q}$, тоді
мають місце порядкові оцінки
\begin{equation}\label{En_teor_pq}
 \mathcal{E}_{Q(N)}(S^{\Omega}_{p,\theta}B)_q\asymp E_{Q(N)}(S^{\Omega}_{p,\theta}B)_q\asymp \frac{1}{N}\big(\log_2 N\big)^{(d-1)\left(\frac{1}{q}-\frac{1}{\theta}\right)_+},
\end{equation}
де $a_+=\max\{0, a\}$.
\rm
\vskip 1 mm
 Зауважимо, що умова  $\Omega\in \Phi_{\alpha,l}$, з деяким $\alpha> \beta= \frac{1}{p}-\frac{1}{q}$ забезпечує те, що для
${f\in S^{\Omega}_{p,\theta}B}$ маємо $f\in S^{\Omega_1}_{q,\theta}B\subset L_q$,
$\Omega_1(\boldsymbol{t})=\Omega(\boldsymbol{t})\boldsymbol{t}^{-\beta}$ і $\|f\|_{S^{\Omega_1}_{q,\theta}B}\ll\|f\|_{S^{\Omega}_{p,\theta}B}$ (див., наприклад,~\cite{Yanchenko_YMG_2010_1}).

\vskip 1mm
\textit{\textbf{Доведення.}}  Спочатку встановимо в
(\ref{En_teor_pq}) оцінки зверху. Нехай $f\in S^{\Omega}_{p,\theta}B$ і ${1<p<q<\infty}$. Тоді, скориставшись співвідношенням
(\ref{e=e}) і лемою~В, можемо записати
$$
E_{Q(N)}(f)_q\asymp
\mathcal{E}_{Q(N)}(f)_q=\bigg\|f(\cdot)-\sum\limits_{\boldsymbol{s}\in
\kappa(N)}\delta_{\boldsymbol{s}}^*(f,\cdot)\bigg\|_q= \bigg\|\sum\limits_{\boldsymbol{s}\in
\kappa^{\perp}(N)}\!\!\!\delta_{\boldsymbol{s}}^*(f,\cdot)\bigg\|_q\!\!\ll
$$
\begin{equation}\label{pqinfty}
\ll\left(\sum \limits_{\boldsymbol{s}\in \kappa^{\perp}(N) }\|\delta_{\boldsymbol{s}}^*
 (f,\cdot)\|^{q}_p 2^{\|\boldsymbol{s}\|_1\left(\frac{1}{p}-\frac{1}{q}\right)q} \right)^{\frac {1}{q}}.
\end{equation}

Далі розглянемо декілька співвідношень між параметрами $q$ i $\theta$.

1) Нехай $q < \theta$. Тоді, для $1 <  \theta < \infty$, застосувавши до $(\ref{pqinfty})$ нерівність Гельдера з показником $\frac{\theta}{q}$ та співвідношення (\ref{1n_log}), одержуємо
$$
 E_{Q(N)}(f)_q\ll\left( \sum \limits_{\boldsymbol{s}\in \kappa^{\perp}(N) }\|\delta_{\boldsymbol{s}}^*
 (f,\cdot)\|^{\theta}_p \big(\Omega(2^{-\boldsymbol{s}})\big)^{-\theta}\right)^{\frac {1}{\theta}}\times
 $$
 $$
 \times\left( \sum \limits_{\boldsymbol{s}\in \kappa^{\perp}(N) }
 \left(2^{\|\boldsymbol{s}\|_1\left(\frac{1}{p}-\frac{1}{q}\right)}\;
 \Omega(2^{-\boldsymbol{s}})\right)^{\frac{q\theta}{\theta-q}}\right)^{\frac {1}{q}-\frac
 {1}{\theta}} \leq
 $$
 $$
 \leq \|f\|_{S_{p,\theta}^{\Omega}B} \left( \sum \limits_{\boldsymbol{s}\in \kappa^{\perp}(N) }
 \left(2^{\|\boldsymbol{s}\|_1\left(\frac{1}{p}-\frac{1}{q}\right)}\;
 \Omega(2^{-\boldsymbol{s}})\right)^{\frac{q\theta}{\theta-q}}\right)^{\frac {1}{q}-\frac
 {1}{\theta}}\ll
 $$
 $$
 \ll \frac{1}{N}|\Theta(N)|^{\frac{1}{q}-\frac{1}{\theta}} \asymp \frac{1}{N}\big(\log_2 N\big)^{(d-1)\left(\frac{1}{q}-\frac{1}{\theta}\right)}
 $$

Якщо ж $\theta=\infty$, то для $f\in S^{\Omega}_{p,\infty}B$,
 згідно з теоремою~А, має місце співвідношення $\|\delta^*_{\boldsymbol{s}}(f,\cdot)\|_p \ll \Omega(2^{-\boldsymbol{s}})$. Тому, скориставшись (\ref{1n_log}),  будемо мати
 $$
 E_{Q(N)}(f)_q\ll \left( \sum \limits_{\boldsymbol{s}\in \kappa^{\perp}(N) }\left(2^{\|\boldsymbol{s}\|_1\left(\frac{1}{p}-\frac{1}{q}\right)}
 \Omega(2^{-\boldsymbol{s}})\right)^q \right)^{\frac {1}{q}} \ll
 \frac{1}{N} |\Theta(N)|^{\frac {1}{q}} \asymp
 $$
 $$
 \asymp \frac{1}{N}\big(\log_2 N\big)^{\frac{d-1}{q}}.
 $$

2)  Нехай тепер $1\leqslant \theta \leqslant q<\infty$, $q\neq 1$. Скориставшись нерівністю
 $$
 \left( \sum \limits_k|a_k|^{v_2}\right)^{\frac{1}{v_2}}\leqslant
 \left( \sum \limits_k|a_k|^{v_1}\right)^{\frac{1}{v_1}}, \
 0<v_1\leqslant v_2<\infty \ ,
 $$
(див.,~\cite[с.~43]{Xardi}), беручи до уваги, що $\Omega \in
\Phi_{\alpha,l}$ з $\alpha>\frac{1}{p}-\frac{1}{q}$ та врахувавши (\ref{Omega_1n}), із (\ref{pqinfty}) отримуємо
 $$
 E_{Q(N)}(f)_q \ll \left( \sum \limits_{\boldsymbol{s}\in \kappa^{\perp}(N)} \!\! \|\delta _{\boldsymbol{s}}^*
 (f,\cdot)\|^{\theta}_p \big(\Omega(2^{-\boldsymbol{s}})\big)^{-\theta}
 2^{\|\boldsymbol{s}\|_1\left(\frac{1}{p}-\frac{1}{q}\right)\theta}
 \big(\Omega(2^{-\boldsymbol{s}})\big)^{\theta}\right)^{\frac {1}{\theta}}\!\!\! \ll
 $$
 $$
 \ll \left( \sum \limits_{\boldsymbol{s}\in \kappa^{\perp}(N)}\|\delta _{\boldsymbol{s}}^*
 (f,\cdot)\|^{\theta}_p \big(\Omega(2^{-\boldsymbol{s}})\big)^{-\theta}
 \right)^{\frac {1}{\theta}} \sup \limits_{\boldsymbol{s}\in \kappa^{\perp}(N)}2^{\|\boldsymbol{s}\|_1\left(\frac{1}{p}-\frac{1}{q}\right)}
 \Omega(2^{-\boldsymbol{s}})\leqslant
 $$
 $$
 \leqslant  \|f\|_{S^{\Omega}_{p,\theta}B}\sup \limits_{\boldsymbol{s}\in \kappa^{\perp}(N)} 2^{\|\boldsymbol{s}\|_1\left(\frac{1}{p}-\frac{1}{q}\right)}
 \Omega(2^{-\boldsymbol{s}}) \ll
 \frac{1}{N}.
 $$

Оцінки зверху в теоремі встановлено.

Перейдемо до встановлення оцінок знизу. Для цього при певних значеннях параметрів $p$, $q$ і $\theta$ достатньо вказати функції $f\in S^{\Omega}_{p,\theta}B$, для яких оцінки
знизу величин $\mathcal{E}_{Q(N)}(f)_q$ співпадають за порядком з
оцінками знизу величин
$\mathcal{E}_{Q(N)}(S^{\Omega}_{p,\theta}B)_q$ в (\ref{En_teor_pq}).
Спочатку означимо функцію, на основі якої буде здійснюватися
побудова таких функцій~$f$.

Для $\boldsymbol{x}=(x_1,\ldots,x_d)$ покладемо
 $$
 D_{\boldsymbol{k}}(\boldsymbol{x})=\prod \limits_{j=1}^{d}D_{k_j}(x_j), \ \
  \boldsymbol{k} \in \mathbb{Z}_+^d,
 $$
 де
 $$
 D_{k_j}(x_j)=\sqrt{\frac {2}{\pi}}\ \left(2\sin {\frac  {x_j}{2}}
 \cos{\frac {2k_j+1}{2}x_j} \right) \cdot x_j^{-1}.
 $$

У роботі \cite{WangHeping_SunYongsheng_1995} показано, що для
перетворення  Фур'є функції $D_{\boldsymbol{k}}(\boldsymbol{x})$
справедлива рівність
 $$
 \mathfrak{F}D_{\boldsymbol{k}}(\boldsymbol{x})=\chi_{\boldsymbol{k}}(\boldsymbol{x})=
 \prod \limits_{j=1}^d \chi_{k_j}(x_j),
 $$
де
$$
      \chi_{n}(\lambda_j)=
 \begin{cases}
    1, & a_j^{n-1}<|\lambda_j|<a_j^{n}, \\
    \frac{1}{2}, & |\lambda_j|=a_j^{n-1} \ \mbox{або} \ |\lambda_j|=a_j^{n}, \\
    0 & \mbox{--- в інших випадках},
 \end{cases}
$$
$$
   \chi_{0}(x_j)=
 \begin{cases}
    1, & |\lambda_j|<1; \\
    \frac{1}{2}, & |\lambda_j|=1; \\
    0, & |\lambda_j|>1.
 \end{cases}
$$

Для оберненого перетворення будемо мати
 $$
 \mathfrak{F}^{-1}\chi_{\boldsymbol{k}}(\boldsymbol{t})=D_{\boldsymbol{k}}(\boldsymbol{x}).
 $$

 Зазначимо, що має місце оцінка \cite{WangHeping_SunYongsheng_1995}
 \begin{equation}\label{Dk}
 \bigg\|\sum \limits_{\boldsymbol{k}\in \rho_+(\boldsymbol{s})}
 D_{\boldsymbol{k}}(\cdot)\bigg\|_p\asymp
 2^{\|\boldsymbol{s}\|_1\left(1-\frac{1}{p}\right)},
 \end{equation}
де
$$
\rho_+(\boldsymbol{s}):=\Big\{\boldsymbol{k}=(k_1,...,k_d): \
\eta(s_j)2^{s_j-1}\leqslant k_j<2^{s_j},  \
k_{j}\in\mathbb{Z}_+^{d},\ j=\overline{1,d}\Big\}.
$$

Далі розглянемо декілька випадків у залежності від значень параметрів $p$, $q$ і $\theta$.

Нехай $\theta =\infty$. Розглянемо функцію
 $$
 f_{1}(\boldsymbol{x})=C_{7}\sum
 \limits_{\boldsymbol{s}\in \Theta(N)}\Omega(2^{-\boldsymbol{s}})2^{-\|\boldsymbol{s}\|_1\left(1-\frac{1}{p}\right)}
 \sum \limits_{\boldsymbol{k} \in \rho_+(\boldsymbol{s})}D_{\boldsymbol{k}}(\boldsymbol{x}).
 $$
  При певному виборі сталої $C_{7}>0$ дана функція належить до класу $S^{\Omega}_{p,\theta}B$ оскільки, скориставшись оцінкою (\ref{Dk}), можемо записати
 $$
 \|f_{1}(\cdot)\|_{S^{\Omega}_{p,\infty}B}=\sup
 \limits_{\boldsymbol{s}\in \Theta(N)}\frac{\|\delta_{\boldsymbol{s}}^*(f_{1},\cdot)\|_p}{\Omega(2^{-\boldsymbol{s}})}=
 $$
 $$
 =\sup
 \limits_{\boldsymbol{s}\in \Theta(N)}C_7 \frac{\Big\|\Omega(2^{-\boldsymbol{s}})2^{-\frac{\|\boldsymbol{s}\|_1}{p'}}
 \sum \limits_{\boldsymbol{k} \in \rho_+(\boldsymbol{s})}D_{\boldsymbol{k}}(\cdot)\Big\|_p}{\Omega(2^{-\boldsymbol{s}})}
 \leqslant
  C_8, \ C_8>0.
 $$

Для $\boldsymbol{s}\in \mathbb{Z}^d_+$ покладемо
 $$
 \Delta(\boldsymbol{s})=\Big\{\boldsymbol{x}: 2^{-s_j-1}\leqslant x_j<2^{-s_j},\ j=\overline {1,d}
 \Big\},
 $$
і зауважимо, що $\Delta(\boldsymbol{s}) \cap \Delta(\boldsymbol{s}')=\varnothing$, якщо $\boldsymbol{s}\neq
\boldsymbol{s}'$. Таким чином, беручи до уваги, що
$S_{Q(N)}(f_{1},\cdot)=0$, скориставшись теоремою~Б та врахувавши, що (див., наприклад, \cite{WangHeping_SunYongsheng_1995})
$$
\bigg|\sum \limits_{k\in\rho(s)}D_k(x)\bigg|=\bigg|\sum
\limits_{k\in\rho(s)}\prod \limits_{j=1}^{d} D_{k_j}(x_j)\bigg| =
\bigg|\prod \limits_{j=1}^{d}\sum
\limits_{k=\eta(s_j)2^{s_j-1}}^{2^{s_j}-1}D_{k_j}(x_j)\bigg|=
$$
$$
=\bigg| \prod \limits_{j=1}^{d}\sqrt{\frac {2}{\pi}} \frac{\sin
2^{s_j}x_j-\sin\eta(s_j)2^{s_j-1}x_j}{x_j}\bigg|,
$$
а також (\ref{Omega_1n}) і (\ref{log}), будемо мати
 $$
  \mathcal{E}_{Q(N)}\big(S^{\Omega}_{p,\theta}B\big)_q\geqslant \mathcal{E}_{Q(N)}(f_{1})_q=
  \|f_{1}(\cdot)\|_q \gg
 $$
 $$
 \gg \Bigg \| \left(
  \sum \limits_{\boldsymbol{s}\in \Theta(N)}|\delta^*_{\boldsymbol{s}}(f_{1},\cdot)|^2
  \right)^{\frac{1}{2}} \Bigg \|_q\geqslant \Bigg(\sum \limits_{\boldsymbol{s}\in \Theta(N)}
  \int \limits_{\Delta(\boldsymbol{s})}|\delta^*_{\boldsymbol{s}}(f_{1},\boldsymbol{x})|^q d\boldsymbol{x}\Bigg)^{\frac{1}{q}}\gg
 $$
 $$
\gg \Bigg(\sum \limits_{\boldsymbol{s}\in \Theta(N)}
  \int \limits_{\Delta(\boldsymbol{s})}\Big|\Omega(2^{-\boldsymbol{s}})2^{-\|\boldsymbol{s}\|_1
  \left(1-\frac{1}{p}\right)}
 \sum \limits_{\boldsymbol{k} \in \rho_+(\boldsymbol{s})}D_{\boldsymbol{k}}(\boldsymbol{x})\Big|^q  d\boldsymbol{x}\Bigg)^{\frac{1}{q}}=
 $$
$$
= \Bigg(\sum \limits_{\boldsymbol{s}\in \Theta(N)}
  \left(\Omega(2^{-\boldsymbol{s}})2^{-\|\boldsymbol{s}\|_1
  \left(1-\frac{1}{p}\right)}\right)^q \int \limits_{\Delta(\boldsymbol{s})}\Big|
 \sum \limits_{\boldsymbol{k} \in \rho_+(\boldsymbol{s})}D_{\boldsymbol{k}}(\boldsymbol{x})\Big|^q  d\boldsymbol{x}\Bigg)^{\frac{1}{q}}\gg
$$
$$
\gg \left( \sum \limits_{\boldsymbol{s}\in \Theta(N)} \left(\Omega(2^{-\boldsymbol{s}})
  2^{\|\boldsymbol{s}\|_1\left(\frac{1}{p}-\frac{1}{q}\right)}
  \right)^q\right)^{\frac{1}{q}}\geqslant \frac{1}{2^l N} |\Theta(N)|^{\frac{1}{q}}\asymp.
 $$
 \begin{equation}\label{ter4f}
\asymp \frac{1}{N} \big(\log_2 N\big)^{\frac{d-1}{q}}.
 \end{equation}

Нехай тепер $1\leqslant\theta\leqslant q< \infty$, $q\neq 1$. Розглянемо функцію
$$
 f_{2}(\boldsymbol{x}):=C_9\Omega(2^{-\tilde{\boldsymbol{s}}})
 2^{-\|\tilde{\boldsymbol{s}}\|_1 \left(1-\frac{1}{p}\right)}\sum \limits_{\boldsymbol{k} \in
  \rho_+(\tilde{\boldsymbol{s}})}D_{{k}}(\boldsymbol{x}), \ \
\tilde{\boldsymbol{s}}\in\Theta(N),\ \ C_9>0.
$$

Згідно з (\ref{Dk}) маємо
$$
\|f_{2}(\cdot)\|_{S^{\Omega}_{p,\theta}B}\asymp\left(\sum\limits_{\boldsymbol{s} \in
  \Theta(N)} \big(\Omega(2^{-\boldsymbol{s}})\big)^{-\theta}\|\delta_{\boldsymbol{s}}^*(f_{2},\cdot)\|_p^{\theta}\right)^{\frac{1}{\theta}}\ll
$$
$$
\ll \left(\big(\Omega(2^{-\tilde{\boldsymbol{s}}})\big)^{-\theta}
\big(\Omega(2^{-\tilde{\boldsymbol{s}}})\big)^{\theta}
 2^{-\theta\|\tilde{\boldsymbol{s}}\|_1\left(1-\frac{1}{p}\right)}\bigg\|\sum \limits_{\boldsymbol{k} \in
  \rho_+(\tilde{\boldsymbol{s}})}D_{\boldsymbol{k}}(\cdot)\bigg\|_p^{\theta}\right)^{\frac{1}{\theta}}\asymp
$$
$$
\asymp  \left( 2^{-\theta\|\tilde{\boldsymbol{s}}\|_1\left(1-\frac{1}{p}\right)}
 2^{\theta\|\tilde{\boldsymbol{s}}\|_1\left(1-\frac{1}{p}\right)}\right)^{\frac{1}{\theta}}= 1,
$$
а отже,  $f_{2}\in S^{\Omega}_{p,\theta}B$ при певному значенні сталої $C_9$.

Врахувавши, що $S_{Q(N)}(f_{2},\cdot)=0$, (\ref{Dk}) та (\ref{Omega_1n}), отримуємо
$$
\mathcal{E}_{Q(N)}(S^{\Omega}_{p,\theta}B)_q\geqslant
\mathcal{E}_{Q(N)}(f_2)_q=\|f_{2}(\cdot)\|_q\gg
$$
$$
\gg \Omega(2^{-\tilde{\boldsymbol{s}}})
 2^{-\|\tilde{\boldsymbol{s}}\|_1 \left(1-\frac{1}{p}\right)}\bigg\|\sum \limits_{\boldsymbol{k} \in
  \rho_+(\tilde{\boldsymbol{s}})}D_{{k}}(\cdot)\bigg\|_q\asymp
$$
$$
\asymp \Omega(2^{-\tilde{\boldsymbol{s}}})
 2^{-\|\tilde{\boldsymbol{s}}\|_1 \left(1-\frac{1}{p}\right)} 2^{\|\tilde{\boldsymbol{s}}\|_1 \left(1-\frac{1}{q}\right)}=\Omega(2^{-\tilde{\boldsymbol{s}}})
 2^{\|\tilde{\boldsymbol{s}}\|_1 \left(\frac{1}{p}-\frac{1}{q}\right)} \geqslant \frac{1}{2^l N}.
$$

У випадку   $1<q<\theta<\infty$ для функції
 $$
  f_{3}(\boldsymbol{x})=C_{10} |\Theta(N)|^{-\frac{1}{\theta}}
  \sum\limits_{\boldsymbol{s} \in
  \Theta(N)}\Omega(2^{-\boldsymbol{s}})
  2^{-\|\boldsymbol{s}\|_1 \left(1-\frac{1}{p}\right)} \sum \limits_{\boldsymbol{k} \in
  \rho_+(\boldsymbol{s})}D_{{k}}(\boldsymbol{x}),
 $$
 скориставшись співвідношенням ($\ref{Dk}$), отримаємо
 $$
  \|f_{3}(\cdot)\|_{S^{\Omega}_{p,\theta}B}\asymp
  |\Theta(N)|^{-\frac{1}{\theta}} \left( \sum \limits_{\boldsymbol{s} \in
  \Theta(N)} \big(\Omega(2^{-\boldsymbol{s}})\big)^{-\theta}
  \|\delta_{\boldsymbol{s}}^*(f_{3},\cdot)\|^{\theta}_{p}
  \right)^{\frac{1}{\theta}} \asymp
 $$
 $$
  \asymp
  |\Theta(N)|^{-\frac{1}{\theta}} \left( \sum \limits_{\boldsymbol{s} \in
  \Theta(N)} \big(\Omega(2^{-\boldsymbol{s}})\big)^{-\theta}\big(\Omega(2^{-\boldsymbol{s}})\big)^{\theta}
  2^{-\|\boldsymbol{s}\|_1 \theta \left(1-\frac{1}{p}\right)} \times\right.
 $$
 $$
 \left. \times \bigg\|\sum \limits_{\boldsymbol{k} \in
  \rho_+(\boldsymbol{s})}D_{{k}}(\cdot)
  \bigg\|^{\theta}_{p}
  \right)^{\frac{1}{\theta}} \ll
 $$
 $$
 \ll |\Theta(N)|^{-\frac{1}{\theta}} \left( \sum \limits_{\boldsymbol{s} \in
  \Theta(N)}
  2^{-\|\boldsymbol{s}\|_1 \theta \left(1-\frac{1}{p}\right)}  2^{\|\boldsymbol{s}\|_1  \theta \left(1-\frac{1}{p}\right)}
  \right)^{\frac{1}{\theta}} = 1.
 $$

Отже, $f_{3}\in S^{\Omega}_{p,\theta}B$ для деякого значення $C_{10}>0$.

Враховуючи, що $S_{Q(N)} (f_{3},\cdot)=0$, та провівши міркування,
аналогічні до тих, що використовувалися для встановлення оцінки
(\ref{ter4f}), одержимо
 $$
 \mathcal{E}_{Q(N)}(S^{\Omega}_{p,\theta}B)_q\geqslant
\mathcal{E}_{Q(N)}(f_{3})_q=\|f_{3}(\cdot)\|_q\gg\frac{1}{N} |\Theta(N)|^{\frac{1}{q}-\frac{1}{\theta}}\asymp
 $$
 $$
\asymp \frac{1}{N}  \big(\log_2 N\big)^{(d-1)\left(\frac{1}{q}-\frac{1}{\theta}\right)}.
 $$

 Оцінки знизу в  (\ref{En_teor_pq}) встановлено. Теорему~1 доведено.

На завершення зробимо коментарі щодо одержаних результатів.

Нехай $\omega(\tau)$~--- функція однієї змінної,
$\omega\in\Phi_{\alpha,l}$, $\alpha>0$, і функція типу мішаного модуля неперервності порядку $l$ задається таким чином
$$
\Omega(\boldsymbol{t})=\Omega(t_1,\ldots,
t_d)=\omega\bigg(\prod\limits_{j=1}^d t_j\bigg), \ \
\Omega\in\Phi_{\alpha,l}, \ \  \alpha>0.
$$
При такому вигляді функціонального параметру $\Omega$  оцінки величини $E_{\bar{Q}_n}\big(S^{\Omega}_{p,\theta}B\big)_q$, у випадку $1<p<q<\infty$, $\alpha>\frac{1}{p}-\frac{1}{q}$, де $\bar{Q}_n=\bigcup\limits_{\|\boldsymbol{s}\|_1<n}Q(\boldsymbol{s})$, встановлено у роботі~\cite{Yanchenko_YMG_2010_1} і, зокрема, у випадку $\Omega(\boldsymbol{t})=\prod\limits_{j=1}^d t_j^{r}$, $\frac{1}{p}-\frac{1}{q}<r<l$~--- в~\cite{WangHeping_SunYongsheng_1995}. Зазначимо, що в~\cite{WangHeping_SunYongsheng_1995} розглядався також випадок, коли $\boldsymbol{r}=(r_1,\ldots,r_d)$, $r_j>\frac{1}{p}-\frac{1}{q}$, $j=\overline{1,d}$.

\vskip 5 mm
Автори висловлюють вдячність А.\,С.~Романюку та В.\,С.~Романюку за їх увагу до роботи та обговорення одержаних результатів.

\vskip 3 mm

\textbf{Contact information:}
Department of the Theory of Functions, Institute of Mathematics of the National
Academy of Sciences of Ukraine, 3, Tereshenkivska st., 01004, Kyiv, Ukraine.

\vskip 3 mm

E-mail: \href{mailto:Yan.Sergiy@gmail.com}{Yan.Sergiy@gmail.com}, \href{mailto:sergstas@ukr.net}{sergstas@ukr.net}


\vskip 3.5mm

\begin{thebibliography}{10}

\bibitem{Stasuk_Yanchenko_Anal_math} S. А. Stasyuk, S.Ya. Yachenko,\emph{Approximation of functions from Nikolskii--Besov type classes of generalized mixed smoothness}~// Anal. Math., \textbf{41}, (2015), 311--334.

\bibitem{Yanchenko_YMG_2016_12} С. Я. Янченко, \emph{Порядкові оцінки апроксимативних характеристик функцій із класів $S^{\Omega}_{p,\theta}B(\mathbb{R}^d)$ із заданою мажорантою мішаних модулів неперервності у рівномірній метриці}~//  Укр. мат. журн., \textbf{68}, (2016), №~12, 1705--1714.

\bibitem{Bari_Stechkin} Н. К. Бари, С. Б. Стечкин, \emph{Наилучшие
приближения и дифференциальные свойства двух сопряженных функций}~//
Тр. Моск. мат. о-ва., \textbf{5}, (1956), 483--522.

\bibitem{Nikolsky_63} С. М. Никольский, \emph{Функции с доминирующей смешанной производной,
 удовлетворяющей кратному условию  Гельдера}~// Сиб. мат. журн., \textbf{4}, (1963), №~6, 1342--1364.

\bibitem{Amanov_1965} Т. И. Аманов, \emph{Теоремы представления и вложения для функциональных пространств
$S^{(r)}_{p,\theta}B(\mathbb{R}_n)$ и $S^{(r)_*}_{p,\theta}B$,
($0\leqslant x_j\leqslant 2\pi$; $j=1,\ldots,n$)}~// Тр. Мат. ин-та
АН СССР., \textbf{77}, (1965), 5--34.

\bibitem{Amanov_book} Т. И. Аманов, \emph{Пространства дифференцируемых функций с доминирующей смешанной производной}, Казах. ССР, Алма-Ата, Наука,  1976.

\bibitem{WangHeping_SunYongsheng_1995}  Wang Heping, Sun Yongsheng,
\emph{Approximation of multivariate functions with a certain mixed
smoothness by entire functions}~// Northeast. Math. J., \textbf{11}, (1995), №~4, 454--466.

\bibitem{WangHeping_2004} Heping Wang,
\emph{Representation and approximation of multivariate function with
bounded mixed smoothness by hyperbolic wavelets}~// J.
Math. Anal. Appl., \textbf{291}, (2004), 698--715.


\bibitem{Temlyakov_1986m} В. Н. Темляков, \emph{Приближение
функций с ограниченной смешанной производной}~// Тр. Мат. ин-та АН
СССР., \textbf{178}, (1986), 1--112.


\bibitem{Romanyuk_2012m} А. С. Романюк, \emph{Аппроксиматиные характеристики классов периодеческих функций многих переменных}, Киев, Ин-тут математики НАН Украины, 2012.


\bibitem{Pustovoitov_94} Н. Н. Пустовойтов, \emph{Представление и приближение периодических
функций многих переменных с заданным смешанным модулем
непрерывности}~// Anal. Math., \textbf{20}, (1994), 35--48.

\bibitem{Pustovoitov_2012} Н. Н. Пустовойтов, \emph{О поперечниках по Колмогорову классов функций с заданным смешанным модулем непрерывности}~// Anal. Math., \textbf{38}, (1994), №~1,  41--64.

\bibitem{WangHeping_SunYongsheng_1997} Wang Sun Yongsheng and Wang Heping,
 \emph{Representation and approximation of multivariate
periodic functions with bounded mixed moduli of smoothness}~// Тр.
Мат. ин-та РАН., \textbf{219}, (1997), 356--377.


\bibitem{Liqin_Duan2010} Liqin Duan, \emph{The best $m$-term approximations on generalized Besov
classes $\boldsymbol{M B^{\Omega}_{q,\theta}}$ with regard to orthogonal dictionaries}~// J. of Approx. Theory, \textbf{162}, (2010), 1964--1981.

\bibitem{Lizorkin_Nikolsky_1989} П. И. Лизоркин, С. М. Никольский, \emph{Пространства функций смешанной гладкости с декомпозиционной точки зрения}~// Тр. Мат. ин-та АН СССР, \textbf{187}, (1989), 143--161.

\bibitem{Lizorkin_69} П. И. Лизоркин,
\emph{Обобщенное лиувиллевское дифференцирование и метод мультипликаторов
в теории вложений классов дифференцируемых функций}~// Тр. Мат. ин-та
АН СССР, \textbf{105}, (1969), 89--167.

\bibitem{D_Zung_86} Динь Зунг, \emph{Приближение функций многих переменных на торе тригонометрическими полиномами}~// Мат. сб., \textbf{131(173)}, (1986), №~2(10), 251--271


\bibitem{Romanyuk_91} А. С. Романюк, \emph{Приближение классов Бесова периодических функций многих переменных в
пространстве $L_q$}~// Укр. мат. журн., \textbf{43}, (1991), №~10, 1398--1408.

\bibitem{Romanyuk_92} А. С. Романюк, \emph{О приближении классов периодических функций многих переменных}~// Укр. мат. журн., \textbf{44}, (1992), №~5, 662--672.

\bibitem{Pustovoitov_99} Н. Н. Пустовойтов,  \emph{Приближение многомерных функций с заданной
мажорантой смешанных модулей непрерывности}~// Мат. заметки, \textbf{65}, (1999), №~1, 107--117.

\bibitem{Stasyuk_2010mz} C. А. Стасюк, \emph{Наилучшие приближения периодических функций многих переменных из классов $B^{\Omega}_{p,\theta}$}~// Мат. заметки, \textbf{87}, (2010), №~1, 108--121.


\bibitem{Stasyuk_2014Tr} C. А. Стасюк, \emph{Приближение суммами Фурье и колмогоровские поперечники классов $\boldsymbol{MB}^{\Omega}_{p,\theta}$ периодических функций нескольких  переменных}~// Тр. Ин-та математики и механики УрО РАН, \textbf{20}, (2014), №~1, 247--257.

\bibitem{Balgimbaeva_2015Tr} Ш. А. Балгимбаева, Т. И. Смирнов, \emph{Оценки поперечников Фурье классов периодических функцй со смешанным модулем гладкосты}~// Тр. Ин-та математики и механики УрО РАН, \textbf{21} , (2015), №~4, 78--94.


\bibitem{Cross_2016} Dinh D\~{u}ng, Vladimir N.~Temlyakov, Tino Ullrich,  \emph{Hyperbolic Cross Approximation}, arXiv:1601.03978v3 [math.NA] 21 Apr 2017.

\bibitem{Nikolsky_1969_book}  С. М. Никольский, \emph{Приближение
функций многих переменных и теоремы вложения}, М.: Наука, 1969.


\bibitem{Lizorkin_1967} П. И. Лизоркин, \emph{Теорема типа Литтльвуда–Палея для кратных интегралов Фурье}~//
Тр. Мат. ин-та АН СССР, \textbf{89}, (1967), 214--230.

\bibitem{Yanchenko_YMG_2010_1} С. Я. Янченко, \emph{Наближення класів
$B^{\Omega}_{p,\theta}$ функцій багатьох змінних у просторі
$L_q(\mathbb{R}^d)$}~// Укр. мат. журн., \textbf{62}, (2010), №~1, 123--135.

\bibitem{Xardi} Г. Г. Харди, Дж. Е. Литтльвуд, Г. Полиa, \emph{Неравенства}, М.: Изд-во иностр. лит., 1948.

\end{thebibliography}
\end{document}